\documentclass{article}

\usepackage{amssymb,amscd,amsmath,amsthm}
\usepackage[matrix, arrow]{xy}
\usepackage{url}

\newtheorem{Thm}{Theorem}
\newtheorem{Prop}[Thm]{Proposition}
\newtheorem{Lem}[Thm]{Lemma}

\newtheorem*{MThm}{Main Theorem}

\newenvironment{Pf}{\textbf{Proof.}\\}{\hspace{\stretch{1}}$\square$}

\DeclareMathOperator{\Stab}{Stab}

\DeclareMathOperator{\Hom}{Hom}

\DeclareMathOperator{\End}{End}
\DeclareMathOperator{\Aut}{Aut}

\newcommand{\bZ}{\mathbb Z}
\newcommand{\bQ}{\mathbb Q}

\title{From Triangulated Categories to Lie Algebras: A Theorem of Peng and Xiao}
\author{Andrew Hubery}
\date{}

\begin{document}

\maketitle

In his seminal article \cite{Ringel}, Ringel showed how to associate to any finitary ring $\Lambda$ an associative unital algebra $\mathcal H(\Lambda)$, with structure constants encoding information about extensions between finite modules. This generalised the Hall algebra \cite{Hall, Steinitz}, which deals with the ring of $p$-adic integers $\bZ_p$ and finite $p$-groups.

In the subsequent article \cite{Ringel2} it is shown that if $\Lambda$ is a representation-directed algebra over a finite field $k$, then the structure constants are given by evaluating integer polynomials. Using these Hall polynomials as structure constants, one may therefore form the generic Ringel-Hall algebra over $\bZ[T]$. Let $\mathfrak n(\Lambda)$ be the subgroup of $\mathcal H(\Lambda)$ generated by the indecomposable modules. If we specialise $T\mapsto1$, then $\bZ\otimes_{\bZ[T]}\mathfrak n(\Lambda)$ becomes a Lie subalgebra of $\bZ\otimes_{\bZ[T]}\mathcal H(\Lambda)$. In fact, over the rational numbers, $\bQ\otimes_{\bZ[T]}\mathcal H(\Lambda)$ is isomorphic to the universal enveloping algebra of $\bQ\otimes_{\bZ[T]}\mathfrak n(\Lambda)$.

In particular, let $\Lambda$ be a representation-finite hereditary $k$-algebra and let $\mathfrak g=\mathfrak n_-\oplus\mathfrak h\oplus\mathfrak n_+$ be the semisimple complex Lie algebra of the same type as $\Lambda$. Then $\bZ\otimes_{\bZ[T]}\mathfrak n(\Lambda)$ can be identified with the Chevalley $\bZ$-form of $\mathfrak n_+$, and $\bZ\otimes_{\bZ[T]}\mathcal H(\Lambda)$ becomes the Kostant $\bZ$-form of the universal enveloping algebra $U(\mathfrak n_+)$ \cite{Ringel3}.

For a general finite dimensional hereditary $k$-algebra $\Lambda$ one considers the composition algebra, the subalgebra generated by the simple modules. This also has a generic version, constructed as a subalgebra of a direct product over infinitely many finite fields of composition algebras \cite{Ringel4}. Green showed in \cite{Green} that the generic composition algebra (after twisting the multiplication via the Euler form of the category $\bmod\Lambda$) is isomorphic to the quantum group of the same type as $\Lambda$.

Therefore, we can realise the quantum group of any symmetrisable Kac-Moody Lie algebra via the module categories of finite dimensional hereditary $k$-algebras.

A natural question is whether it is possible to obtain the full (quantised) enveloping algebra, or at least the full Lie algebra. The latter question was answered by Peng and Xiao in \cite{PX1} for the affine Lie algebras of type $\widetilde{\mathbb A}$, and in \cite{PX2} for the simple complex Lie algebras, using the root category $\mathcal D^b(\bmod\Lambda)/T^2$. Finally, these methods were generalised to all 2-periodic triangulated $k$-categories (satisfying some finiteness conditions), and in particular to the root category of a finite dimensional hereditary $k$-algebra \cite{PX3}.

In particular, given any symmetrisable generalised Cartan matrix, the (derived) Kac-Moody Lie algebra can be realised via the root categories (together with the Grothendieck groups) of finite dimensional hereditary $k$-algebras.

Unfortunately, it is unlikely that the approach taken by Peng and Xiao yields an associative algebra.

In this article, we offer a simplified and more intuitive proof of the theorem of Peng and Xiao \cite{PX3}. We begin by providing a more categorical proof of the associativity of multiplication for the standard Ringel-Hall algebras. Until now the proof has relied on counting filtrations, and so cannot easily be adapted to triangulated categories. The new proof exhibited in Section 2 works in the setting of an exact category, with the key result being Proposition \ref{P1}. Replacing the pull-back/push-out diagram by the Octahedral Axiom immediately yields the analogous Proposition \ref{P2} for triangulated categories. This clarifies and improves the results in Section 6 of \cite{PX3}.

In Section 4 we offer a unified approach to the calculations in Section 7 of \cite{PX3}, and then use these results to prove the Jacobi identity. We note that in Case (I), an extra argument is required to prove that coefficient of $\tilde h_X$ vanishes. This argument is missing in \cite{PX3}.

In the final section we consider when we can endow the Lie algebra with a symmetric bilinear form. Contrary to \cite{PX3}, this is not always possible, even for the root category of a finite dimensional hereditary $k$-algebra, since it is not in general defined on the Cartan subalgebra. The best we can hope for is to define the form on $\mathfrak h\times\mathfrak h_1$, where $\mathfrak h_1$ is a subgroup of the Cartan subalgebra. The proof that this is invariant then follows from our previous considerations on the Jacobi identity.

We remark that To\"en has recently shown how to construct an associative algebra from a dg-category \cite{Toen} (under some finiteness assumptions), and it will clearly be of interest to investigate the connections between his derived Hall algebra and the quantised enveloping algebras.

{\bf Acknowledgements} These notes were prepared for a series of lectures at Universit\"at Paderborn. I would like to thank all involved for helpful comments and suggestions, and in particular H.~Asashiba, W.~Crawley-Boevey, B.~Deng, H.~Krause and J.~Xiao.

\section{The Main Theorem}

Let $k$ be a finite field and let $\mathcal T$ be an idempotent complete triangulated $k$-category. We shall also assume that $\mathcal T$ is skeletally small, so the isomorphism classes of objects form a set, and 2-periodic, i.e. $T^2\cong1$, where $T$ is the shift of $\mathcal T$. Note that for $X$ indecomposable, $\End X$ is a finite dimensional local $k$-algebra. Set $d(X):=\dim(\End X/\mathrm{rad}\End X)$.

The Grothendieck group $\mathcal G$ of $\mathcal T$ is the quotient of the free abelian group with generators the isomorphism classes of objects in $\mathcal T$ by the relations $[X]+[Y]-[L]$ whenever there exists an exact triangle $Y\to L\to X\to TY$. Denote by $h_X$ the image of $[X]$ in $\mathcal G$. Note that $h_{TM}=-h_M$ for all objects $M$, and that the group is generated by all $h_X$ for $X$ indecomposable. The Grothendieck group is called proper if $h_X\neq0$ in $\mathfrak h$ for all indecomposable objects $X$.

We call $\mathcal T$ finitary if all homomorphism spaces are finite sets and, given a fixed element $h\in\mathcal G$, there are only finitely many isomorphism classes of indecomposable objects $X$ (up to shift) with $h_X=h$.

\begin{MThm}
Let $k$ be a finite field with $q_k$ elements and let $\mathcal T$ be a skeletally small and idempotent complete triangulated $k$-category. Assume also that $\mathcal T$ is finitary and 2-periodic with proper Grothendieck group. Then we can associate to $\mathcal T$ a Lie algebra $\mathfrak g(\mathcal T)$ over the ring $\bZ/(q_k-1)\bZ$.
\end{MThm}

The rest of the article is devoted to proving this theorem.

\section{The Ringel-Hall Algebra of an Exact Category}

We recall the definition of an exact category. Let $\mathcal A$ be an additive category and let $\mathcal E$ be a class of kernel-cokernel pairs $(f,g)$, closed under isomorphism. We call $f$ an inflation, $g$ a deflation and $(f,g)$ a conflation.

The pair $(\mathcal A,\mathcal E)$ is an exact category in the sense of Quillen \cite{Quillen} (see also the appendix to \cite{Keller}) if the following axioms hold .
\begin{enumerate}
\item[(Ex0)] $0\xrightarrow{1} 0$ is a deflation.
\item[(Ex1)] The composition of two deflations is again a deflation.
\item[(Ex2)] For $f:Y\to L$ and a deflation $m:M\to L$ there exists a pull-back
$$\xymatrix{L'\ar@{.>}[r]^{f'}\ar@{.>}[d]^{m'} &M\ar[d]^m\\Y\ar[r]^f &L}$$
with $m'$ a deflation.
\item[(Ex3)] For $m':L'\to Y$ and an inflation $f':L'\to M$ there exists a push-out
$$\xymatrix{L'\ar[r]^{f'}\ar[d]^{m'} & M\ar@{.>}[d]^m\\Y\ar@{.>}[r]^f &L}$$
with $f$ an inflation.
\end{enumerate}

We recall the following result from \cite{Keller}.
\begin{Lem}
The product of two inflations is again an inflation. Moreover, the pair $((f'\, -m')^t,(m\,f))$ obtained from either (Ex2) or (Ex3) is a conflation.
\end{Lem}

The Grothendieck group $\mathcal G$ of $\mathcal A$ is the quotient of the free abelian group with generators the isomorphism classes of objects of $\mathcal A$ by the relations $[X]+[Y]-[L]$ whenever there exists a conflation $Y\to L\to X$. Denote by $h_X$ the image of $[X]$ in $\mathcal G$.

We call $\mathcal A$ finitary if each homomorphism group is finite and, given any $h\in\mathcal G$, there exist only finitely many isomorphism classes $[X]$ with $h_X=h$. N.B. The category of finite modules over a finitary ring, in the sense of Ringel \cite{Ringel2}, is then a finitary abelian category.

Let $\mathcal A$ be a finitary and skeletally small exact category and let $W_{XY}^L$ denote the set of all conflations $Y\to L\to X$. The group $\Aut(X,Y):=\Aut X\times\Aut Y$ acts on $W_{XY}^L$ via
$$\xymatrix{Y\ar[r]^f\ar[d]^\eta &L\ar[r]^g\ar@{=}[d] &X\ar[d]^\xi\\
Y\ar[r]^{\overline f} &L\ar[r]^{\overline g} &X}$$
and we denote the quotient set by $V_{XY}^L$. Since $f$ is monic and $g$ epic this action is free, so
$$F_{XY}^L:=\left|V_{XY}^L\right|=\frac{\left|W_{XY}^L\right|}{\left|\Aut(X,Y)\right|}.$$

The Ringel-Hall algebra $\mathcal H(\mathcal A)$ is defined as follows. Form the free $\bZ$-module on generators indexed by the set of isomorphism classes of objects, writing $u_X$ for $u_{[X]}$, and use the numbers $F_{XY}^L$ as structure constants. That is,
$$u_Xu_Y:=\sum_{[L]}F_{XY}^Lu_L.$$
We note that this sum is finite and that $u_0$ is a unit for the multiplication. We now prove that the multiplication is associative.

Let us define an action of $\Aut(X,Y,Z,L):=\Aut X\times\Aut Y\times\Aut Z\times\Aut L$ on the pairs of conflations $W_{XY}^L\times W_{LZ}^M$ via
$$\xymatrix{Y\ar[r]^f\ar[d]^\eta &L\ar[r]^g\ar[d]^\lambda &X\ar[d]^\xi\\
Y\ar[r]^{\overline f} &L\ar[r]^{\overline g} &X} \qquad
\xymatrix{Z\ar[r]^l\ar[d]^\zeta &M\ar[r]^m\ar@{=}[d] &L\ar[d]^\lambda\\
Z\ar[r]^{\overline l} &M\ar[r]^{\overline m} &L}$$
Dually we have an action of $\Aut(X,Y,Z,L')$ on the pairs $W_{XL'}^M\times W_{YZ}^{L'}$.

We fix a complete set of representatives $\mathcal M$ for all isomorphism classes of objects of $\mathcal A$.

\begin{Prop}\label{P1}
For $X,Y,Z,M\in\mathcal M$ there is a bijection between the sets
$$\bigcup_{L\in\mathcal M}\frac{W_{XY}^L\times W_{LZ}^M}{\Aut(X,Y,Z,L)}\longleftrightarrow\bigcup_{L'\in\mathcal M}\frac{W_{XL'}^M\times W_{YZ}^{L'}}{\Aut(X,Y,Z,L')}.
$$
\end{Prop}

\begin{Pf}
Consider the pair $((f,g),(l,m))\in W_{XY}^L\times W_{LZ}^M$. Let $L'\in\mathcal M$ be the representative of the pull-back
$$\xymatrix{L'\ar[r]^{f'}\ar[d]^{m'}&M\ar[d]^m\\Y\ar[r]^f&L}$$
Take $l':Z\to L'$ such that $f'l'=l$ and $m'l'=0$. Then $l'$ is a kernel for the deflation $m'$, so $(l',m')$ is a conflation. Similarly, setting $g':=gm$, then $g'$ is a deflation and $f'$ is a kernel for $g'$, so $(f',g')$ is a conflation. Thus we have a commutative diagram
$$\xymatrix{Z\ar[d]^{l'}\ar@{=}[r] &Z\ar[d]^l\\L'\ar[d]^{m'}\ar[r]^{f'} &M\ar[d]^m\ar[r]^{g'} &X\ar@{=}[d]\\Y\ar[r]^f &L\ar[r]^g &X}$$
with all rows and columns being conflations. In particular, we obtain the pair $((f',g'),(l',m'))$ in $W_{XL'}^M\times W_{YZ}^{L'}$.

Now consider the pair $((\overline f,\overline g),(\overline l,\overline m))=(\xi,\eta,\zeta,\lambda)\cdot((f,g),(l,m))$. Using the same construction as above, we obtain $((\overline f',\overline g'),(\overline l',\overline m'))$ in $W_{XL''}^M\times W_{YZ}^{L''}$, where $L''\in\mathcal M$ is the representative of the pull-back of $\overline m$ along $\overline f$.

There exists $\lambda'$ giving a morphism of conflations
$$\xymatrix@C=15pt{L'\ar[rr]^-{(f'\,-m')^t}\ar@{.>}[d]^{\lambda'} &&M\oplus Y\ar[rr]^-{(m\,f)}\ar[d]^-{\left(\begin{smallmatrix}1\\&\eta\end{smallmatrix}\right)} &&L\ar[d]^\lambda\\
L''\ar[rr]^-{(\overline f'\,-\overline m')^t} &&M\oplus Y\ar[rr]^-{(\overline m\,\overline f)} &&L}$$
It follows that $\lambda'$ is an isomorphism, hence $L''=L'$ and $\lambda'\in\Aut L'$.

Similarly we obtain $\xi'$ and $\zeta'$ giving isomorphisms of conflations
$$\xymatrix{L'\ar[r]^{f'}\ar[d]^{\lambda'} &M\ar[r]^{g'}\ar@{=}[d] &X\ar@{.>}[d]^{\xi'}\\
L'\ar[r]^{\overline f'} &M\ar[r]^{\overline g'} &X} \qquad
\xymatrix{Z\ar[r]^{l'}\ar@{.>}[d]^{\zeta'} &L'\ar[r]^{m'}\ar[d]^{\lambda'} &Y\ar[d]^\eta\\
Z\ar[r]^{\overline l'} &L'\ar[r]^{\overline m'} &Y}$$
Thus $((\overline f',\overline g'),(\overline l',\overline m'))=(\xi',\eta,\zeta',\lambda')\cdot((f',g'),(l',m'))$ and so the two pairs of conflations lie in the same $\Aut(X,Y,Z,L')$-orbit.

In fact, since $l'$ is monic and $g'$ epic, we deduce that $\zeta'=\zeta$ and $\xi'=\xi$.

This proves that the map from left to right is well-defined. Using the symmetry of the situation, we obtain the required bijection.
\end{Pf}

Let $N_{XYZ}^{LM}$ denote the size of the quotient set $\frac{W_{XY}^L\times W_{LZ}^M}{\Aut(X,Y,Z,L)}$. Since the action is again free, we have
$$N_{XYZ}^{LM}=\frac{\left|W_{XY}^L\times W_{LZ}^M\right|}{\left|\Aut(X,Y,Z,L)\right|}=F_{XY}^L F_{LZ}^M.$$
Hence we have the associativity of the multiplication
$$\sum_{[L]}F_{XY}^LF_{LZ}^M=\sum_{[L']}F_{XL'}^MF_{YZ}^{L'}.$$

\begin{Thm}
The Ringel-Hall algebra $\mathcal H(\mathcal A)$ of a finitary and skeletally small exact category $\mathcal A$ is an associative, unital algebra.
\end{Thm}

Suppose further that $\mathcal A$ is a $k$-category for some finite field $k$. Set $q_k:=\left|k\right|$ and write $\mathfrak n$ for the subgroup of $\mathcal H(\mathcal A)$ generated by the $u_X$ for $X$ indecomposable.

\begin{Thm}
Over the ring $\bZ/(q_k-1)$, $\mathfrak n$ becomes a Lie subalgebra of $\mathcal H(\mathcal A)$.
\end{Thm}

\begin{Pf}
It is enough to show that $F_{XY}^L-F_{YX}^L\equiv 0\mod (q_k-1)$ whenever $X$ and $Y$ are indecomposable and $L$ is decomposable (c.f. \cite{Ringel4}). We may assume that $X\not\cong Y$, since otherwise the result is trivial. So suppose that $L=L'\oplus L''$ and consider a conflation $Y\xrightarrow{(f'\,f'')^t} L'\oplus L''\xrightarrow{(g'\,g'')} X$.

If $L\not\cong X\oplus Y$ then $k^*$ acts freely on the set $V_{XY}^L$, hence $F_{XY}^L\equiv0$. For, let $\big(\eta,\big(\begin{smallmatrix}1\\&\kappa\end{smallmatrix}\big),\xi\big)\in\Aut(Y,L,X)$ with $\kappa\in k^*$ be an automorphism of this conflation. Since $\mathcal A$ is idempotent complete, $\End Y$ is a finite dimensional local $k$-algebra. As neither $f'$ nor $f''$ is 0, both $\eta-1$ and $\eta-\kappa$ are nilpotent. Therefore $\kappa-1$ is also nilpotent, so $\kappa=1$.

Now suppose that $L=X\oplus Y$, so that every conflation is split. We may assume that $f''=1$, $g'=1$ and $f''=-g'=\theta:Y\to X$. Thus $V_{XY}^L\cong\Hom(Y,X)$ and hence $F_{XY}^{X\oplus Y}\equiv 1\mod (q_k-1)$. Similarly $F_{YX}^{X\oplus Y}\equiv1$ and we are done.
\end{Pf}

\section{Triangulated Categories}

We would like to imitate this result for a triangulated category $\mathcal T$, and in particular for the root category of a finite dimensional hereditary algebra.

We first recall the definition of a triangulated category. Let $\mathcal T$ be an additive category endowed with an equivalence $T$  and let $\mathcal E$ be a class of triangles $Y\xrightarrow{f} L\xrightarrow{g} X\xrightarrow{h} TL$ closed under isomorphism. We call the triangles in $\mathcal E$ exact.

The pair $(\mathcal T,\mathcal E)$ is a triangulated category if the following axioms hold \cite{Puppe,Verdier} (see also \cite{Krause,Neeman}).
\begin{enumerate}
\item[(Tr0)] The triangle $0\to X\xrightarrow{1} X\to 0$ is exact for all $X$.
\item[(Tr1)] Each map $f$ fits into an exact triangle $(f,g,h)$.
\item[(Tr2)] The triangle $(f,g,h)$ is exact if and only if $(g,h,-Tf)$ is exact.
\item[(Tr3)] Given two exact triangles $(f,g,h)$ and $(f',g',h')$ and two maps $\eta$ and $\lambda$ such that $\lambda f=f'\eta$, there exists a morphism $(\eta,\lambda,\xi)$ of exact triangles
$$\xymatrix{Y\ar[r]^f\ar[d]^\eta &L\ar[r]^g\ar[d]^\lambda &X\ar[r]^h\ar@{.>}[d]^\xi &TY\ar[d]^{T\eta}\\
Y'\ar[r]^{f'} &L'\ar[r]^{g'} &X'\ar[r]^{h'} &TY'}$$
\item[(Tr$4'$)] For $f:Y\to L$ and $m:M\to L$ there exists a morphism of triangles
$$\xymatrix{L'\ar[r]^{f'}\ar[d]^{m'} &M\ar[r]^{g'}\ar[d]^m &X\ar[r]^{h'}\ar@{=}[d] &TL'\ar[d]^{Tm'}\\
Y\ar[r]^f &L\ar[r]^g &X\ar[r]^h &TY}$$
such that the left hand square is homotopy cartesian with differential $\delta=h'g$. That is, the triangle $((f'\,-m')^t,(m\,f),h'g)$ is exact.
\end{enumerate}

We shall need the following results about triangulated categories. (In fact, none of the following needs the axiom (Tr$4'$)). A triangle split if it is isomorphic to a triangle of the form $Y\xrightarrow{(0\,1)^t} X\oplus Y\xrightarrow{(1\,0)} X\xrightarrow{0} TY$. 

\begin{Lem}\label{LemA}
Let $Y\xrightarrow{f} L\xrightarrow{g} X\xrightarrow{h} TY$ be a triangle.
\begin{enumerate}
\item If in (Tr3) $\eta$ and $\lambda$ are isomorphisms, then so is $\xi$;
\item This triangle is split if and only if $h=0$;
\item Suppose $L=M\oplus P$ and that $f=(f'\,0)^t$. Then $(f,g,h)$ is isomorphic to $Y\xrightarrow{(f'\,0)^t} M\oplus P\xrightarrow{\big(\begin{smallmatrix}g'\\&1\end{smallmatrix}\big)} N\oplus P\xrightarrow{(h'\,0)} TY$ and $(f',g',h')$ is exact. A similar result also holds whenever $g=(g'\,0)$.
\end{enumerate}
\end{Lem}

\begin{Pf}
The first part is Proposition 1.1.20 in \cite{Neeman}, and the second part is Corollary 1.2.7.

The third occurs as Lemma 2.5 in \cite{PX2}, but the proof can be simplified as follows. Form an exact triangle $Y\xrightarrow{f'} M\xrightarrow{g'} N\xrightarrow{h'} TY$. The direct sum of triangles $Y\xrightarrow{(f'\,0)^t} M\oplus P\xrightarrow{\big(\begin{smallmatrix}g'\\&1\end{smallmatrix}\big)} N\oplus P\xrightarrow{(h'\,0)} TY$ is exact by Proposition 1.2.1 of \cite{Neeman}, and there exists a morphism $(1,1,\xi)$ from $(f,g,h)$ to this triangle, which must be an isomorphism.
\end{Pf}

We recall that (Tr4') is equivalent to the following axiom (Tr4), or the Octahedral Axiom (see \cite{Neeman} or the appendix in \cite{Krause}).
\begin{enumerate}
\item[(Tr4)] Given exact triangles $(f,g,h),(l,m,n)$ and $(l',m',n')$ with $n'=nf$, there exists an exact triangle $(f',g',h')$ such that the following diagram commutes
$$\xymatrix{Z\ar@{=}[r]\ar[d]^{l'} &Z\ar[d]^l\\
L'\ar@{.>}[r]^{f'}\ar[d]^{m'} &M\ar@{.>}[r]^{g'}\ar[d]^m &X\ar@{.>}[r]^-{h'}\ar@{=}[d] &TL'\ar[d]^{Tm'}\\
Y\ar[r]^f\ar[d]^{n'} &L\ar[r]^g\ar[d]^n &X\ar[r]^-{h} &TY\\
TZ\ar@{=}[r] &TZ}$$
Moreover, the middle left square is homotopy cartesian with differential $\delta=h'g=(Tl')n$.
\end{enumerate}

We call $\mathcal T$ finitary if each homomorphism group is finite and, given any element $h$ of the Grothendieck group $\mathcal G$ of $\mathcal T$, there exist only finitely many isomorphism classes of indecomposables $[X]$ (up to shift) with $h_X=h$. In this case we have another characterisation of split triangles: the triangle $(f,g,h)$ is split if and only if $L\cong X\oplus Y$. For, we can apply the cohomological functor $\Hom(-,TY)$. By counting arguments, $\Hom(h,TY)=0$, so $h=0$.

Now let $\mathcal T$ be a finitary and skeletally small triangulated category. We prove that the analogue of Proposition \ref{P1} still holds for $\mathcal T$.

Denote by $W_{XY}^L$ the set of all exact triangles $Y\to L\to X\to TY$. The group $\Aut(X,Y)=\Aut X\times \Aut Y$ acts via
$$\xymatrix{Y\ar[d]^\eta\ar[r]^f &L\ar@{=}[d]\ar[r]^g &X\ar[d]^\xi\ar[r]^-h &TY\ar[d]^{T\eta}\\
Y\ar[r]^{\overline f} &L\ar[r]^{\overline g} &X\ar[r]^-{\overline h} &TY}$$
and we write $V_{XY}^L$ for the quotient. Note that this action is not free.

Let $\Aut(X,Y,Z,L)$ act on $W_{XY}^L\times W_{LZ}^M$ via
$$\xymatrix{Y\ar[r]^f\ar[d]^\eta &L\ar[r]^g\ar[d]^\lambda &X\ar[r]^-h\ar[d]^\xi &TY\ar[d]^{T\eta}\\
Y\ar[r]^{\overline f} &L\ar[r]^{\overline g} &X\ar[r]^-{\overline h} &TY} \qquad
\xymatrix{Z\ar[r]^l\ar[d]^\zeta &M\ar[r]^m\ar@{=}[d] &L\ar[r]^-n\ar[d]^\lambda &TZ\ar[d]^{T\zeta}\\
Z\ar[r]^{\overline l} &M\ar[r]^{\overline m} &L\ar[r]^-{\overline n} &TZ}$$
Dually, we have an action of $\Aut(X,Y,Z,L')$ on the set $W_{XL'}^M\times W_{YZ}^{L'}$.

We fix a complete set of representatives $\mathcal M$ for all isomorphism classes of objects of $\mathcal T$.

\begin{Prop}\label{P2}
For $X,Y,Z,M\in\mathcal M$ there is a bijection between the sets
$$\bigcup_{L\in\mathcal M}\frac{W_{XY}^L\times W_{LZ}^M}{\Aut(X,Y,Z,L)}\longleftrightarrow\bigcup_{L'\in\mathcal M}\frac{W_{XL'}^M\times W_{YZ}^{L'}}{\Aut(X,Y,Z,L')}.$$
\end{Prop}

\begin{Pf}
Consider the pair $((f,g,h),(l,m,n))\in W_{XY}^L\times W_{LZ}^M$. Applying (Tr4), we obtain for some $L'\in\mathcal M$ a pair $((f',g',h'),(l',m',n'))\in W_{XL'}^M\times W_{YZ}^{L'}$ fitting into a commutative diagram such that the middle left square is homotopy cartesian with differential $\delta=h'g=(Tl')n$.

Now consider the pair $((\overline f,\overline g,\overline h),(\overline l,\overline m,\overline n))=(\xi,\eta,\zeta,\lambda)\cdot((f,g,h),(l,m,n))$. Using (TR4) we again find for some $L''\in\mathcal M$ a pair $((\overline f',\overline g',\overline h'),(\overline l',\overline m',\overline n'))\in W_{XL''}^M\times W_{YZ}^{L''}$ fitting into a commutative diagram such that the middle left square is homotopy cartesian with differential $\overline\delta=\overline h'\overline g=(T\overline l')\overline n$.

Using (Tr3) we can find $\lambda'$ giving a morphism of exact triangles
$$\xymatrix@C=15pt{L'\ar[rr]^-{(f'\,-m')^t}\ar@{.>}[d]^{\lambda'} &&M\oplus Y\ar[rr]^-{(m\,f)}\ar[d]^-{\left(\begin{smallmatrix}1\\&\eta\end{smallmatrix}\right)} &&L\ar[rr]^\delta\ar[d]^\lambda &&TL'\ar@{.>}[d]^{T\lambda'}\\
L''\ar[rr]^-{(\overline f'\,-\overline m')^t} &&M\oplus Y\ar[rr]^-{(\overline m\,\overline f)} &&L\ar[rr]^{\overline\delta} &&TL''}$$
It follows that $\lambda'$ is an isomorphism, hence $L''=L'$ and $\lambda'\in\Aut L'$.

Similarly we obtain $\xi'$ and $\zeta'$ giving isomorphisms of exact triangles
$$\xymatrix{L'\ar[d]^{\lambda'}\ar[r]^{f'} &M\ar@{=}[d]\ar[r]^{g'} &X\ar@{.>}[d]^{\xi'}\ar[r]^{h'} &TL'\ar[d]^{T\lambda'}\\L'\ar[r]^{\overline f'} &M\ar[r]^{\overline g'} &X\ar[r]^{\overline h'} &TL'} \qquad
\xymatrix{Z\ar@{.>}[d]^{\zeta'}\ar[r]^{l'} &L'\ar[d]^{\lambda'}\ar[r]^{m'} &Y\ar[d]^\eta\ar[r]^{n'} &TZ\ar@{.>}[d]^{T\zeta'}\\Z\ar[r]^{\overline l'} &L'\ar[r]^{\overline m'} &Y\ar[r]^{\overline n'} &TZ}$$

We deduce that $((\overline f',\overline g',\overline h'),(\overline l',\overline m',\overline n'))=(\xi',\eta,\zeta',\lambda')\cdot((f',g',h'),(l',m',n'))$ and so the two pairs of exact triangles lie in the same $\Aut(X,Y,Z,L')$-orbit.

This shows that the map from left to right is well-defined. Using the symmetry of the situation, we obtain the required bijection.
\end{Pf}

\textbf{Remark.}
In fact, it follows from the proof that $\bar h'(\xi'-\xi)g=0$ and $\bar l'(\zeta'-\zeta)T^{-1}n=0$. However, unlike for an exact category, we cannot deduce that $\xi'=\xi$ and $\zeta'=\zeta$. Similarly, none of the actions described above is free. In particular, it is not clear whether we can define an associative algebra using the numbers $F_{XY}^L:=\left|V_{XY}^L\right|$ as structure constants.

\section{Reduction modulo $q_k-1$}

From now on, let $k$ be a finite field with $q_k$ elements, and $\mathcal T$ a finitary and skeletally small triangulated $k$-category. Furthermore, we assume that $\mathcal T$ has split idempotents. Thus the endomorphism ring of an indecomposable object is a finite dimensional local $k$-algebra. For an indecomposable $X$, denote by $d(X)$ the dimension $\dim_k(\End X/\mathrm{rad}\End X)$.

\textbf{Remark.}
This is not a serious restriction on our triangulated category, since we can always form the idempotent completion \cite{BS}.

For fixed indecomposable objects $X,Y$ and $Z$ we will calculate the numbers $F_{XY}^L, F_{LZ}^M$ and $N_{XYZ}^{LM}:=\left|\frac{W_{XY}^L\times W_{LZ}^M}{\Aut(X,Y,Z,L)}\right|$. It will be convenient to work over the subring $R\subset\bQ$ obtained by localising $\bZ$ at the set of numbers coprime to $(q_k-1)$.

The method of proof here and later is based on the following observation. If a finite group $G$ acts on a finite set $X$, then
$$\frac{\left|X\right|}{\left|G\right|}=\sum_{x\in X/G}\frac{1}{\left|\Stab x\right|}.$$

\subsection{The numbers $F_{XY}^L$}

\begin{Lem}\label{L1}
\begin{enumerate}
\item If $L\not\cong M\oplus TZ$, then $F_{LZ}^M\equiv \frac{\left|W_{LZ}^M\right|}{\left|\Aut(Z,L)\right|} \mod (q_k-1)$.
\item If $L\cong M\oplus TZ$, then $F_{LZ}^M=1$ and $\left|W_{LZ}^M\right|=\frac{\left|\Aut L\right|}{\left|\Hom(TZ,M)\right|}$.
\end{enumerate}
\end{Lem}

\begin{Pf}
Suppose first that $L\not\cong M\oplus TZ$. Using the observation, we shall prove that the stabilisers are all bijective to vector spaces and hence have size 1 modulo $(q_k-1)$.

For $t=(l,m,n)\in W_{LZ}^M$ define
$$S(t):=\{(\zeta-1,0,\lambda-1)\in\End t \mid (\zeta,\lambda)\in\Stab(t)\}.$$
Here we have written $\End t$ for the subalgebra of $\End Z\times\End M\times\End L$ consisting of all triples $(\zeta,\mu,\lambda)$ giving a morphism of the triangle $t$.

Clearly $S(t)$ is in bijection with $\Stab(t)$ and contains 0. For $(\bar\zeta,0,\bar\lambda)\in S(t)$ we know that $\bar\zeta$ is either nilpotent or an isomorphism, and that $l\bar\zeta=0$. Since $l\neq 0$, $\bar\zeta$ is nilpotent.

Let $(\bar\zeta_i,0,\bar\lambda_i)$ for $i=1,2$ be two elements of $S(t)$ and let $\alpha\in k$. Then $\bar\zeta:=\bar\zeta_1+\alpha\bar\zeta_2$ is nilpotent, and hence $\zeta:=1+\bar\zeta$ is an isomorphism. Set $\lambda:=1+\bar\lambda_1+\alpha\bar\lambda_2$, so that $(\zeta,1,\lambda)\in\End t$. Thus $\lambda$ is an isomorphism and so $(\zeta,\lambda)\in\Stab(t)$. Therefore $S(t)$ is a vector subspace of $\End t$. This proves Part 1.

Now suppose that $L=M\oplus TZ$. There is a unique orbit, represented by the triangle $Z\xrightarrow{0} M\xrightarrow{(1\,0)^t} M\oplus TZ\xrightarrow{(0\,1)} TZ$. Clearly $(\zeta,\lambda)$ stablises this triangle if and only if $\lambda=\left(\begin{smallmatrix}1&\theta\\0&T\zeta\end{smallmatrix}\right)$ for some $\theta:TZ\to M$.
\end{Pf}

\begin{Lem}\label{L2}
\begin{enumerate}
\item If $L\not\not\cong 0$, then $F_{XY}^L\equiv \frac{\left|W_{XY}^L\right|}{\left|\Aut(X,Y)\right|} \mod (q_k-1)$.
\item If $L\cong X\oplus Y$ and $X\not\cong Y$, then $F_{XY}^L=\frac{\left|W_{XY}^L\right|}{\left|\Aut(X,Y)\right|}=\left|\Hom(Y,X)\right|$.
\item We have $\left|W_{XY}^{X\oplus Y}\right|=\frac{\left|\Aut(X\oplus Y)\right|}{\left|\Hom(X,Y)\right|}$.
\end{enumerate}
\end{Lem}

\begin{Pf}
This first part follows from the previous lemma, since $X\cong L\oplus TY$ and $X$ and $Y$ indecomposable implies $L\cong0$.

Suppose $L=X\oplus Y$ and $X\not\cong Y$. Since every triangle $(f,g,h)$ is split, each orbit is represented by a triangle of the form $((\theta\,1)^t,(1\,-\theta),0)$ for some $\theta:Y\to X$. These orbits are all distinct, and have trivial stabilisers.

Finally consider $W_{XY}^{X\oplus Y}$. By rotation, this set is in bijection with $W_{X\oplus YT^{-1}X}^Y$, whose size we know by the previous lemma.
\end{Pf}

\subsection{The numbers $N_{XYZ}^{LM}$}

\begin{Lem}\label{L3}
If $L\not\cong M\oplus TZ$, then $N_{XYZ}^{LM} \equiv F_{XY}^LF_{LZ}^M$.
\end{Lem}

\begin{Pf}
If $L=0$, then it is clear that $N_{XYZ}^{LM}=F_{XY}^LF_{LZ}^M$, and this equals 0 unless both $X\cong TY$ and $Z\cong M$, in which case it equals 1.

Suppose $L\not\cong 0$ and let $(t,u):=((f,g,h),(l,m,n))\in W_{XY}^L\times W_{LZ}^M$. For $\beta$ an automorphism write $\bar\beta:=\beta-1$ and define
$$S(t,u):=\{((\bar\eta,\bar\lambda,\bar\xi),(\bar\zeta,0,\bar\lambda))\in\End t\times\End u \mid (\xi,\eta,\zeta,\lambda)\in\Stab(t,u)\}.$$
As before, $S(t,u)$ is in bijection with $\Stab(t,u)$ and contains the zero morphism.

For $((\bar\eta,\bar\lambda,\bar\xi),(\bar\zeta,0,\bar\lambda))\in S(t,u)$ we know that $\bar\zeta$ is nilpotent, say $\bar\zeta^{r-1}=0$. Then $n\bar\lambda^{r-1}=0$, so $\bar\lambda^{r-1}=m\theta$ for some $\theta:L\to M$. Hence $\bar\lambda^r=0$, since $\bar\lambda m=0$.

We now have $(\bar\eta^r,0,\bar\xi^r)\in\End t$. If $\bar\xi^r$ is an automorphism, then $g=0$ and the rotation $(-T^{-1}h,f,g)$ of $t$ is split. Since $X$ and $Y$ are indecomposable, $L\cong0$, a contradiction. Hence $\bar\xi$ is nilpotent. Similarly $\bar\eta$ is nilpotent.

Let $((\bar\eta_i,\bar\lambda_i,\bar\xi_i),(\bar\zeta_i,0,\bar\lambda_i))$ for $i=1,2$ be two elements of $S(t,u)$ and let $\alpha\in k$. Set $\bar\xi:=\bar\xi_1+\alpha\bar\xi_2$ and $\xi:=1+\bar\xi$, and similarly for $\eta,\zeta$ and $\lambda$. Then $\bar\xi,\bar\eta$ and $\bar\zeta$ are nilpotent, so $\xi,\eta$ and $\zeta$ are isomorphisms. Since $((\eta,\lambda,\xi),(\zeta,1,\lambda))\in\End t\times\End u$, $\lambda$ must also be an isomorphism and so $(\xi,\eta,\zeta,\lambda)\in\Stab(t,u)$. Therefore $S(t,u)$ is a vector subspace of $\End t\times\End u$.

It follows that $N_{XYZ}^{LM}\equiv\frac{\left|W_{XY}^L\times W_{LZ}^M\right|}{\left|\Aut(X,Y,Z,L)\right|}$, which we know by the previous subsection is congruent to the product $F_{XY}^LF_{LZ}^M$.
\end{Pf}

\begin{Lem}\label{L4}
Let $L=M\oplus TZ$ with $M\not\cong 0$, so that $F_{LZ}^M=1$ and $\left|W_{LZ}^M\right|=\frac{\left|\Aut L\right|}{\left|\Hom(TZ,M)\right|}$.
\begin{enumerate}
\item If $L\not\cong X\oplus Y$, then $N_{XYZ}^{LM} \equiv \frac{\left|W_{XY}^L\right|}{\left|\Aut(X,Y,Z)\right|}$.
\item If $L\cong X\oplus Y$, then there are three subcases.
\begin{itemize}
\item[i)] If $M\not\cong Y$, then $N_{XYZ}^{LM}=1$.
\item[ii)] If $M\not\cong X$, then $N_{XYZ}^{LM}-1\equiv\frac{\dim\Hom(Y,X)}{d(X)}$.
\item[iii)] If $M\cong X\cong Y$, then $N_{XYZ}^{LM}-2 \equiv \frac{\dim\mathrm{rad}\End X}{d(X)}$.
\end{itemize}
\end{enumerate}
\end{Lem}

\begin{Pf}
Every orbit is represented by a pair $(t,u)=\big((f,g,h),(0,(1\,0)^t,(0\,1))\big)$. We shall use the terminology of the previous lemma. Then $\lambda=\left(\begin{smallmatrix}1&\theta\\0&T\zeta\end{smallmatrix}\right)$ for some $\theta:TZ\to M$.

Write $f=(f_1\,f_2)^t$ and $g=(g_1\,g_2)$. Then
$$f_1\bar\eta=\theta f_2, \quad f_2\bar\eta=(T\bar\zeta)f_2, \quad \bar\xi g_1=0, \quad \bar\xi g_2=g_1\theta+g_2(T\bar\zeta), \quad \bar\eta h=h\bar\xi.$$

1. If $L\not\cong X\oplus Y$ then $\bar\xi,\bar\eta$ and $\bar\zeta$ are all nilpotent. For, if $\bar\xi$ is an automorphism, then $g_1=0$ and $Y\cong M\oplus Y'$ by Lemma \ref{LemA}. Since $Y$ is indecomposable and $M\not\cong 0$, we must have that $Y\cong M$ and $X\cong TZ$, hence $L\cong X\oplus Y$, a contradiction. Thus $\bar\xi$ is nilpotent. Now, if $\bar\eta$ is an automorphism, then $h=0$. Thus $t$ is split and $L\cong X\oplus Y$, a contradiction. Finally, if $\bar\zeta$ is an automorphism, then $f_2=0$ and Lemma \ref{LemA} applies to give $X\cong M$ and $Y\cong TZ$, hence $L\cong X\oplus Y$, a contradiction.

Therefore $S(t,u)$ is a vector subspace of $\End t\times\End u$ and the result follows, using the formula for $\left|W_{LZ}^M\right|$.

2. Now suppose that $L\cong X\oplus Y$. Then $t$ is split and $h=0$. Also, either $f_1$ or $f_2$ is an isomorphism.

Suppose that $f_2$ is an isomorphism, so we may identify $Y=TZ$ and $X=M$. Then there is a unique orbit for $(t,u)$, given by the pair of split triangles. It follows that $\xi=1$, $\theta=0$ and $\eta=T\zeta$, so the stabiliser is isomorphic to $\Aut Y$.

Conversely, suppose that $f_2$ is not an isomorphism. Then we may identify $Y=M$ and $X=TZ$, and $t$ is represented by $((1\,\phi)^t,(-\phi\,1),0)$ with $\phi:Y\to X$ not an isomorphism. We note the conditions
$$\xi\phi=\phi, \quad T\zeta=\xi+\phi\theta, \quad \eta=1+\theta\phi, \quad \lambda=\left(\begin{smallmatrix}1&\theta\\0&T\zeta\end{smallmatrix}\right).$$
Since $\phi$ is not an isomorphism, $\eta$ is necessarily invertible, and $\zeta$ and $\lambda$ are invertible if and only if $\xi$ is. Therefore the stabiliser is determined by $\theta$ and $\xi$ such that $\bar\xi\phi=0$.

If $\phi=0$, then the stabiliser is isomorphic to $\Aut X\times\Hom(X,Y)$. Otherwise, if $\phi\neq0$, then $\bar\xi$ must be nilpotent and the stabiliser is in bijection with some subspace of $\mathrm{rad}\End X\times\Hom(X,Y)$.

We conclude that if $M\not\cong Y$, then $N_{XYZ}^{LM}=1$.

If $M\not\cong X$ then $N_{XYZ}^{LM}-1$ is congruent to
$$\frac{\left|W_{XY}^L\times W_{LZ}^M\right|}{\left|\Aut(X,Y,Z,L)\right|}-\frac{1}{\left|\Aut X\right|\left|\Hom(X,Y)\right|}=\frac{\left|\Hom(Y,X)\right|-1}{\left|\Hom(X,Y)\right|\left|\Aut X\right|}$$
using Lemma \ref{L2}. Since $\End X$ is local, we can write $\left|\Aut X\right|=q^r(q^{d(X)}-1)$, where $r=\dim\mathrm{rad}\End X$. Thus
$$N_{XYZ}^{LM}-1 \equiv \frac{1}{q^{r+\dim\Hom(X,Y)}}\cdot\frac{q^{\dim\Hom(Y,X)}-1}{q^{d(X)}-1} \equiv \frac{\dim\Hom(Y,X)}{d(X)}.$$

Finally, if $M=X=Y=TZ$, then $N_{XYZ}^{LM}-2$ is congruent to
$$\frac{\left|W_{XX}^{X^2}\right|}{\left|\Aut X\right|^3\left|\End X\right|}-\frac{1}{\left|\Aut X\right|\left|\End X\right|}-\frac{1}{\left|\Aut X\right|}.$$
From Lemma \ref{L2}, $\left|W_{XX}^{X^2}\right|=\left|\Aut X^2\right|/\left|\End X\right|$. Using the notation above, $\left|\Aut X^2\right|=q^{4r+d(X)}(q^{d(X)}+1)(q^{d(X)}-1)^2$, so that
$$N_{XYZ}^{LM}-2 \equiv \frac{1}{q^{2r+d(X)}}\cdot\frac{q^r-1}{q^{d(X)}-1}\equiv \frac{r}{d(X)}.$$
\end{Pf}

We state the dual result for the numbers $\hat N_{XYZ}^{ML'}:=\left|\smash[t]{\frac{W_{XL'}^M\times W_{YZ}^{L'}}{\Aut(X,Y,Z,L')}}\right|$, with the special case being when $L'\cong M\oplus T^{-1}X$.

\begin{Lem}\label{L5}
If $L'\not\cong M\oplus T^{-1}X$, then $\hat N_{XYZ}^{ML'}\equiv F_{XL'}^MF_{YZ}^{L'}$.

On the other hand, let $L'=M\oplus T^{-1}X$ with $M\not\cong0$.
\begin{enumerate}
\item If $L'\not\cong Y\oplus Z$, then $\hat N_{XYZ}^{ML'}\equiv\frac{\big|W_{YZ}^{L'}\big|}{\left|\Aut(X,Y,Z)\right|}$.
\item If $L'\cong Y\oplus Z$, then there are three subcases.
\begin{itemize}
\item[i)] If $M\not\cong Y$, then $\hat N_{XYZ}^{ML'}=1$.
\item[ii)] If $M\not\cong Z$, then $\hat N_{XYZ}^{ML'}-1\equiv\frac{\dim\Hom(Z,Y)}{d(Z)}$.
\item[iii)] If $M\cong Y\cong Z$, Then $\hat N_{XYZ}^{ML'}-2\equiv\frac{\dim\mathrm{rad}\End Z}{d(Z)}$.
\end{itemize}
\end{enumerate}
\end{Lem}

\section{The Jacobi Identity}

In order to prove the Jacobi identity, we shall need two further assumptions on our triangulated category $\mathcal T$. Firstly, we need that $\mathcal T$ is 2-periodic, i.e. $T^2\cong1$. We note that this is satisfied for the root category $\mathcal D^b(\bmod\Lambda)/T^2$ for $\Lambda$ a finite dimensional hereditary $k$-algebra \cite{PX2} (see also \cite{Keller2}). We shall also need that for $X$ indecomposable, $h_X\neq0$ in $\mathcal G$. The latter condition is referred to as proper in \cite{PX3}.

We continue to work over the ring $R/(q_k-1)$.

For $X$ indecomposable, define $\tilde h_X:=\frac{h_X}{d(X)}$ in $\bQ\otimes_\bZ\mathcal G$ and set $\mathfrak h$ to be the integer lattice generated by the $\tilde h_X$. Define a bilinear form on $\mathfrak h$ by setting $(h_X\,|\,h_Y)$ to be the alternating sum
$$\dim\Hom(X,Y)-\dim\Hom(X,TY)+\dim\Hom(Y,X)-\dim\Hom(Y,TX).$$
Note that $(\tilde h_X\,|\,h_Y)\in\bZ$ for all indecomposable objects $X$ and $Y$.

Define $\mathfrak n$ to be the free abelian group on generators indexed by the isomorphism classes of indecomposable objects, writing $u_X$ for $u_{[X]}$.

In this section, $X,Y$ and $Z$ will again denote indecomposable objects.

We define a bracket product on $\mathfrak g:=\mathfrak h\oplus\mathfrak n$ via the formulae
\begin{align*}{}
[u_X,u_Y] &:= \begin{cases}
\sum_{L\,\mathrm{ind}}(F_{XY}^L-F_{YX}^L)u_L &\textrm{if }X\not\cong TY;\\
\tilde h_X &\textrm{if }X\cong TY,
\end{cases}\\{}
[\tilde h_X,u_Y] &:= -(\tilde h_X\mid h_Y)u_Y, \quad [u_Y,\tilde h_X]=(\tilde h_X\mid h_Y)u_Y\\{}
[\tilde h_X,\tilde h_Y] &:= 0.
\end{align*}
Clearly this is antisymmetric. It is well defined from the finiteness condition on the Grothedieck group.

We wish to prove the Jacobi identity
$$[[u_X,u_Y],u_Z]-[[u_X,u_Z],u_Y]-[[u_Z,u_Y],u_X]\equiv 0$$
in the following three cases.
\begin{enumerate}
\item[(I)] $X,Y\not\cong TZ$ and $X\not\cong TY$. In this case,
$$[[u_X,u_Y],u_Z]=\sum_{\substack{L,M\,\mathrm{ind}\\L\not\cong TZ}}(F_{XY}^L-F_{YX}^L)(F_{LZ}^M-F_{ZL}^M)u_M - (F_{XY}^{TZ}-F_{YX}^{TZ})\tilde h_Z.$$
\item[(II)] $X,Y\not\cong TZ$ and $X\cong TY$. In this case, $[[u_X,u_Y],u_Z]=-(\tilde h_X\,|\,h_Z)u_Z$.
\item[(III)] $X\cong Y\cong TZ$. In this case, the Jacobi identity is trivial.
\end{enumerate}
N.B. It is enough to consider these cases. For, if $X\cong TZ$ and $Y\not\cong TZ$, then either $X\not\cong TY$ and we can use the triple $(u_Z,u_X,u_Y)$, or $X\cong TY$ and we can use the triple $(u_Y,u_Z,u_X)$. Similarly if $X\not\cong TZ$ and $Y\cong TZ$.

\begin{Lem}\label{L6}
If $L$ is decomposable or 0, then $F_{XY}^L-F_{YX}^L\equiv 0$.
\end{Lem}

\begin{Pf}
For $L=0$ this is clear, so suppose that $L$ is decomposable, say $L=M\oplus TZ$ with $Z$ indecomposable and $M\not\cong 0$. If $L\not\cong X\oplus Y$, then by Lemmas \ref{L2} and \ref{L4},
$$F_{XY}^L\equiv\frac{\left|W_{XY}^L\right|}{\left|\Aut(X,Y)\right|}\equiv\left|\Aut Z\right|N_{XYZ}^{LM}\equiv 0.$$
If $L=X\oplus Y$ but $X\not\cong Y$, then by Lemma \ref{L2}, $F_{XY}^L=\left|\Hom(Y,X)\right|\equiv 1$. Finally, if $X\cong Y$, then the left hand side is clearly 0.
\end{Pf}

This implies that in the expression for $[[u_X,u_Y],u_Z]$ in Case (I), we can remove the restriction that $L$ is indecomposable. Moreover, since the Grothedieck group is proper, we can remove the restriction that $L\not\cong TZ$. For, $F_{TZZ}^M\neq0$ implies $h_M=0$, so that $M$ is decomposable.

\textbf{Case (I)}
The coefficient $c_M$ of $u_M$ in the Jacobi identity equals the alternating sum
$$c_M=\Delta_{XYZ}^M+\Delta_{YZX}^M+\Delta_{ZXY}^M-\Delta_{YXZ}^M-\Delta_{ZYX}^M-\Delta_{XZY}^M,$$
where
$$\Delta_{XYZ}^M:=\sum_LF_{XY}^LF_{LZ}^M-\sum_LF_{XL}^MF_{YZ}^L.$$
N.B. This summation makes sense in $R/(q_k-1)$, since if $L$ is decomposable and not isomorphic to $X\oplus Y$, then $F_{XY}^L\equiv 0$.

Using Lemmas \ref{L3} and \ref{L5}, Proposition \ref{P2} and that $T^2=1$ we see that
$$\Delta_{XYZ}^M \equiv F_{XY}^{M\oplus TZ}F_{M\oplus TZZ}^M-F_{XM\oplus TX}^MF_{YZ}^{M\oplus TX}-N_{XYZ}^{M\oplus TZM}+\hat N_{XYZ}^{MM\oplus TX}.$$
Substituting this into $c_M$, the part involving the $F$s can be expressed as
$$(F_{XY}^{M\oplus TZ}-F_{YX}^{M\oplus TZ})F_{M\oplus TZ Z}^M - F_{XM\oplus TX}^M(F_{YZ}^{M\oplus TX}-F_{ZY}^{M\oplus TX})$$
together with two other terms given by cyclically permuting $(X,Y,Z)$. These all vanish by Lemma \ref{L6}.

It remains to consider the terms involving $N$ and $\hat N$. By our conditions on $(X,Y,Z)$, we must always be in Case 1 of Lemmas \ref{L4} and \ref{L5}. We see that $\frac{\left|W_{XY}^{M\oplus TZ}\right|}{\left|\Aut(X,Y,Z)\right|}$ occurs once from $N_{XYZ}^{M\oplus TZ M}$ and once from $\hat N_{ZXY}^{MM\oplus TZ}$, and these cancel in $c_M$. Similarly for each of the other permutations.

We deduce that the coefficient of $u_M$ in the Jacobi identity is 0.

It remains to consider the terms lying in $\mathfrak h$. That is, we need to check that
$$-(F_{XY}^{TZ}-F_{YX}^{TZ})\tilde h_Z + (F_{XZ}^{TY}-F_{ZX}^{TY})\tilde h_Y + (F_{ZY}^{TX}-F_{YZ}^{TX})\tilde h_X \equiv 0.$$

\begin{Lem}\label{L7}
Let $t:Y\to Z\to X\to TY$ be a triangle. Then $\End t$ is local.
\end{Lem}

\begin{Pf}
We will show that if one of $\xi,\eta$ or $\zeta$ is nilpotent, then they all are.

Suppose that $\xi$ is nilpotent, say $\xi^m=0$. Then $(\eta^m,\zeta^m,0)\in\End t$, and since $t$ is not split, $\eta^m$ must be nilpotent. Similarly, after rotating the triangle, we see that $\zeta$ is nilpotent.
\end{Pf}

Set $d(t):=\dim(\End t/\mathrm{rad}\End t)$. The above proof shows that we have a natural monomorphism
$$\End t/\mathrm{rad}\End t\to\End X/\mathrm{rad}\End X,$$
and hence $d(t)$ divides $d(X)$, and similarly $d(Y)$ and $d(Z)$.

The group $\Aut(X,Y,Z)$ acts on the set $W_{XY}^{TZ}$, so
$$\frac{\left|W_{XY}^{TZ}\right|}{\left|\Aut(X,Y,Z)\right|}=\sum_t\frac{1}{\left|\Aut t\right|},$$
where the sum is taken over all orbits. Applying Lemma \ref{L1}, we get
$$F_{XY}^{TZ}\equiv \frac{\left|W_{XY}^{TZ}\right|}{\left|\Aut(X,Y)\right|}=\frac{\left|W_{XY}^{TZ}\right|\left|\Aut(Z)\right|}{\left|\Aut(X,Y,Z)\right|} = \sum_t\frac{\left|\Aut Z\right|}{\left|\Aut t\right|}\equiv \sum_t\frac{d(Z)}{d(t)}.$$

If $t\in W_{XY}^{TZ}$, then $h_X+h_Y+h_Z=0$ and so
$$\frac{d(X)}{d(t)}\tilde h_X+\frac{d(Y)}{d(t)}\tilde h_Y+\frac{d(Z)}{d(t)}\tilde h_Z=0$$
By rotating triangles, we deduce that
$$-F_{XY}^{TZ}\tilde h_Z \equiv -\sum_t\frac{d(Z)}{d(t)}\tilde h_Z \equiv \sum_t\left(\frac{d(X)}{d(t)}\tilde h_X+\frac{d(Y)}{d(t)}\tilde h_Y\right) \equiv F_{YZ}^{TX}\tilde h_X+F_{ZX}^{TY}\tilde h_Y.$$

This completes the proof of the Jacobi identity in Case (I).

\textbf{Case (II)}
Since the Grothendieck group is proper, $u_X$ does not occur in $[u_X,u_Z]$. Furthermore, all numbers of the form $F_{XTX}^M$ vanish. Therefore
$$-[[u_X,u_Z],u_{TX}]-[[u_Z,u_{TX}],u_X]=\sum_Mc_Mu_M,$$
where
$$c_M=\Delta_{XTXZ}^M+\Delta_{TXZX}^M+\Delta_{ZXTX}^M-\Delta_{TXXZ}^M-\Delta_{ZTXX}^M-\Delta_{XZTX}^M.$$
In particular, this lies entirely in $\mathfrak n$.

As before, the part of $c_M$ involving the $F$s vanishes. For the terms involving $N$ and $\hat N$, we have two cases, depending on whether $M\cong Z$ or not. If $M\not\cong Z$, then we are always in Case 1 of Lemmas \ref{L4} and \ref{L5}, and $c_M\equiv 0$ as in Case (I). So suppose that $M=Z$. Then, for example,
$$N_{XTXZ}^{M\oplus TZ M}\equiv \tfrac{\left|W_{XTX}^{M\oplus TZ}\right|}{\left|\Aut(X,Y,Z)\right|}, \quad N_{TXZX}^{M\oplus TXM}\equiv 1+\tfrac{\dim\Hom(Z,TX)}{d(X)}, \quad N_{ZXTX}^{M\oplus XM}\equiv 1.$$
A short calculation (noting $\Hom(TX,Z)\cong\Hom(X,TZ)$) reveals that $c_Z=(\tilde h_X\,|\,h_Z)$. Since $[[u_X,u_{TX}],u_Z]=-(\tilde h_X \mid h_Z)u_Z$, we are done.

We have therefore shown that the Jacobi identity is satisfied for all triples $(u_X,u_Y,u_Z)$.

Now consider the triple $(u_X,u_Y,\tilde h_Z)$. Then $[[u_X,u_Y],\tilde h_Z]$ equals 0 if $X\cong TY$; otherwise we have
$$[[u_X,u_Y],\tilde h_Z] = \sum_{L\,\mathrm{ind}}(F_{XY}^L-F_{YX}^L)[u_L,\tilde h_Z] = (\tilde h_Z\mid h_X+h_Y)[u_X,u_Y].$$
The Jacobi identity follows immediately. The remaining case $(u_X,\tilde h_Y,\tilde h_Z)$ is easily checked.

This completes the proof of the Main Theorem.

\section{The Symmetric Bilinear Form}

It would clearly be of interest to extend the symmetric bilinear form to an invariant form on the whole of the Lie algebra $\mathfrak g/(q_k-1)\mathfrak g$. Unfortunately, this is not always possible since although $(\tilde h_X\,|\,h_Y)\in\bZ$, it is not true that $(\tilde h_X\,|\,\tilde h_Y)\in\bZ$. In fact, this even fails for the root category $\mathcal D^b(\bmod \Lambda)/T^2$ of a finite dimensional hereditary $k$-algebra $\Lambda$ (for example, when $\Lambda$ is the $\mathbb F_4$-species of type $\mathbb G_2$).

One important special case when we can define the bilinear form is when $\Lambda=kQ$ is the path algebra of a quiver. For, the Grothendieck group is in this case freely generated by the simple modules, and these all have trivial endomorphism rings.

In general, define $\mathfrak h_1$ to be the sublattice of $\mathfrak h$ generated by $h_X$ for $X$ indecomposable. Set $\mathfrak g_1:=\mathfrak h_1\oplus\mathfrak n$ and extend the bilinear form on $\mathfrak h\times\mathfrak h_1$ to the whole of $\mathfrak g\times\mathfrak g_1$ via
$$(\tilde h_X\mid u_Y):=0, \quad (u_X\mid u_Y):=\begin{cases}1 &\textrm{if }X\cong TY;\\0 &\textrm{otherwise}.\end{cases}$$
This form is then invariant over $R/(q_k-1)$ in the sense that
\begin{align*}
(\tilde h_X\mid[u_Y,u_Z])d(Z) &\equiv -([\tilde h_X,u_Y]\mid u_Z)\\
\textrm{and}\quad d(X)([u_X,u_Y]\mid u_Z) &\equiv (u_X\mid[u_Y,u_Z])d(Z).
\end{align*}
For the first condition, both sides are 0 unless $Y\cong TZ$, in which case they are both congruent to $(\tilde h_X|h_Y)$. For the second, we need to show that
$$d(X)\big(F_{XY}^{TZ}-F_{YX}^{TZ}\big)\equiv d(Z)\big(F_{YZ}^{TX}-F_{ZY}^{TX}\big).$$
As in Case (I) of the Jacobi identity, we know that
$$d(X)F_{XY}^{TZ}\equiv\sum_t\frac{d(X)d(Z)}{d(t)}\equiv d(Z)F_{YZ}^{TX},$$
and we are done.

{\small
Andrew Hubery\\
Universit\"at Paderborn, Germany\\
hubery@math.upb.de}


\begin{thebibliography}{99}
\bibitem{BS} P.~Balmer and M.~Schlichting, `Idempotent completion of triangulated categories', \textit{J. Algebra} 236 (2001) 819--834.
\bibitem{Green} J.A.~Green, `Hall algebras, hereditary algebras and quantum groups', \textit{Invent. Math.} 120 (1995) 361--377.
\bibitem{Hall} P.~Hall, `The algebra of partitions', in \textit{Proc. 4th Canad. Math. Congress} (Banff, 1957) (Univ. of Toronto Press, Toronto, 1959) 147--159.
\bibitem{Keller} B.~Keller, `Chain complexes and stable categories', \textit{Manus. Math.} 67 (1990) 379--417.
\bibitem{Keller2} B.~Keller, `On triangulated orbit categories', preprint \url{http://www.math.jussieu.fr/~keller/publ/triaorbit.ps} .
\bibitem{Krause} H.~Krause, `Derived categories, resolutions, and Brown representability', \url{http://wwwmath.upb.de/~hkrause/publ/chicago.ps} .
\bibitem{Quillen} D.~Quillen, \textit{Higher algebraic K-theory, I}, Springer Lecture Notes in Mathematics 341 (1973) 85--147.
\bibitem{Neeman} A.~Neeman, \textit{Triangulated Categories}, Annals of Mathematical Studies 148 (Princeton University Press, Princeton, 2001).
\bibitem{PX1} L.~Peng and J.~Xiao, `A realisation of affine Lie algebras of type $\tilde A_{n-1}$ via the derived categories of cyclic quivers', in \textit{Representation theory of algebras} (Cocoyoc, 1994) CMS Conf. Proc. 18 (Amer. Math. Soc., Providence, 1996) 539--554.
\bibitem{PX2} L.~Peng and J.~Xiao, `Root categories and simple Lie algebras', \textit{J. Algebra} 198 (1997) 19--56.
\bibitem{PX3} L.~Peng and J.~Xiao, `Triangulated categories and Kac-Moody Lie algebras', \textit{Invent. Math.} 140 (2000) 563--603.
\bibitem{Puppe} D.~Puppe, `On the formal structure of stable homotopy theory', in \textit{Colloquium on algebraic topology}, Matematisk Institut, Aarhus Universitet (Aarhus, 1962) 65--71.
\bibitem{Ringel} C.M.~Ringel, `Hall algebras and quantum groups', \textit{Invent. Math.} 101 (1990) 583--591.
\bibitem{Ringel2} C.M.~Ringel, `Hall algebras', in \textit{Topics in Algebra, Part 1}, Banach Center Publ. 26, Part 1 (PWN, Warsaw, 1990) 433--447.
\bibitem{Ringel3} C.M.~Ringel, `Hall polynomials for the representation-finite hereditary algebras', \textit{Adv. Math.} 84 (1990) 137--178.
\bibitem{Ringel4} C.M.~Ringel, `From representations of quivers via Hall and Loewy algebras to quantum groups', \textit{Contemp. Math.} 131 (1992) 381--401.
\bibitem{Steinitz} E.~Steinitz, `Zur Theorie der Abel'schen Gruppen', Jahresber. Deutsch. Math.-Verein 9 (1901) 80--85. 
\bibitem{Toen} B.~To\"en, `Derived Hall algebras', preprint \url{http://front.math.ucdavis.edu/math.QA/0501343} .
\bibitem{Verdier} J.L.~Verdier, `Des cat\'egories d\'eriv\'ees des cat\'egories ab\'eliennes', \textit{Ast\'erique} 239 (1996)
\end{thebibliography}
\end{document}